\documentclass[12pt]{amsart}

\usepackage{epsf}
\usepackage{amsmath,amsthm,amssymb,amsfonts}
\usepackage{colordvi}
\usepackage{color}
\usepackage[all]{xy}

\newcommand{\xyinc}{\ar@{^{(}->}}

\numberwithin{equation}{section}
\numberwithin{figure}{section}
\theoremstyle{plain}
\newtheorem{thm}[equation]{Theorem}
\newtheorem{prop}[equation]{Proposition}
\newtheorem{lem}[equation]{Lemma}
\newtheorem{cor}[equation]{Corollary}

\theoremstyle{definition}
\newtheorem{defi}[equation]{Definition}
\newtheorem{exas}[equation]{Examples} 

\newtheorem{rem}[equation]{Remark} 

\newcounter{FNC}[page]
\def\newfootnote#1{{\addtocounter{FNC}{2}$^\fnsymbol{FNC}$%
     \let\thefootnote\relax\footnotetext{$^\fnsymbol{FNC}$#1}}}

 \newcommand{\ch}{\mathrm{char}}
\newcommand{\sumsub}[1]{\sum_{\substack{#1}}} 

\newcommand{\ten}{\mbox{\raisebox{1pt}{${\scriptstyle \otimes}$}}}
\newcommand{\iten}{\backslash}
\newcommand{\tof}{{\scriptstyle\,\#_f\,}}

\newcommand{\gr}{\mathrm{gr}}
\newcommand{\Sh}{\mathrm{Sh}}     
\newcommand{\Nod}{\mathrm{Nod}}     
\newcommand{\id}{\mathrm{id}}
\newcommand{\End}{\mathrm{End}}

\newcommand{\tsso}{\mathrm{t.s.s.o.}}

\newcommand{\onto}{\twoheadrightarrow} 
\newcommand{\map}[1]{\xrightarrow{#1}}

\newcommand{\QSym}{\mathcal{Q}\mathit{Sym}}
\newcommand{\SSym}{\mathfrak{S}\mathit{Sym}}
\newcommand{\YSym}{\mathcal{Y}\mathit{Sym}}

\newcommand{\field}{\Bbbk}

\newcommand{\calH}{\mathcal{H}}

\newcommand{\GL}{\calH_O}
\newcommand{\HGL}{\calH_{HO}}
\newcommand{\euler}{\mathbf{e}}

\newcommand{\tzero}{\begin{picture}(5,10)(-1,0)\put(1.5,2.5){\circle*{3}}\end{picture}}
\newcommand{\tone}{
\begin{picture}(5,10)(-1,0)\thicklines
\put(1.5,-2){\line(0,1){10}}
           \put(1.5,-2){\circle*{3}}\put(1.5,8){\circle*{3}}\end{picture}}
 \newcommand{\tonebig}{
\begin{picture}(5,12)(-1,0)\thicklines
\put(1.5,-2){\line(0,1){12}}
           \put(1.5,-2){\circle*{3}}\put(1.5,10){\circle*{3}}\end{picture}}
          
\newcommand{\ttwotwo}{
\begin{picture}(20,20)(-1,0)\thicklines
 \put(10,-2){\line(-1,2){10}}\put(0,18){\circle*{3}}
 \put(10,-2){\line(1,2){10}}\put(20,18){\circle*{3}}
 \put(10,-2){\circle*{3}}
       \end{picture}}
\newcommand{\ttwoone}{\begin{picture}(5,20)(-1,0)\thicklines
\put(1.5,-3){\line(0,1){20}}
           \put(1.5,-2){\circle*{3}}\put(1.5,8){\circle*{3}}\put(1.5,18){\circle*{3}}
       \end{picture}}
\newcommand{\tthreefive}{ \begin{picture}(20,20)(-1,0)\thicklines
    \put(10,-2){\line(-1,2){10}} \put( 0,18){\circle*{3}}
    \put(10,-2){\line( 0,1){20}} \put(10,18){\circle*{3}}
    \put(10,-2){\line( 1,2){10}} \put(20,18){\circle*{3}} \put(10,-2){\circle*{3}}
   \end{picture}}
\newcommand{\tthreethree}{\begin{picture}(20,24)(-1,0)\thicklines
\put(9,-2){\line(-3,4){9}}
\put(9,-2){\line(3,4){9}}\put(18,10){\line(0,1){12}}
\put(9,-2){\circle*{3}}\put(18,22){\circle*{3}}\put(0,10){\circle*{3}}
\put(18,10){\circle*{3}}
       \end{picture}}
\newcommand{\tthreefour}{\begin{picture}(20,24)(-1,0)\thicklines
\put(9,-2){\line(-3,4){9}}
\put(9,-2){\line(3,4){9}}\put(0,10){\line(0,1){12}}
\put(9,-2){\circle*{2}}\put(0,22){\circle*{3}}\put(0,10){\circle*{3}}
\put(18,10){\circle*{3}}
       \end{picture}}
\newcommand{\tthreetwo}{\begin{picture}(20,24)(-1,0)\thicklines
\put(9,10){\line(-3,4){9}}
\put(9,10){\line(3,4){9}}\put(9,-2){\line(0,1){12}}
\put(9,-2){\circle*{3}}\put(9,10){\circle*{3}}\put(0,22){\circle*{3}}
\put(18,22){\circle*{3}}
       \end{picture}}
\newcommand{\tthreeone}{\begin{picture}(5,36)(-1,0)\thicklines
\put(1.5,-2){\line(0,1){36}}
\put(1.5,-2){\circle*{3}}\put(1.5,10){\circle*{3}}
\put(1.5,22){\circle*{3}}\put(1.5,34){\circle*{3}}
\end{picture}}
\newcommand{\tfourone}{\begin{picture}(21,24)(-2,0)\thicklines
\put(9,-2){\line(-3,4){9}}
\put(9,-2){\line(3,4){9}}\put(9,-2){\line(0,1){12}}
\put(18,10){\line(0,1){12}}\put(18,22){\circle*{3}}\put(9,-2){\circle*{3}}
\put(9,10){\circle*{3}}\put(0,10){\circle*{3}}\put(18,10){\circle*{3}}
       \end{picture}}

\newcommand{\tfourthree}{\begin{picture}(28,24)(0,0)\thicklines
       \put(18,-2){\line(-3,4){9}}\put(18,-2){\line(3,4){9}}
\put(9,10){\line(-3,4){9}}\put(9,10){\line(3,4){9}}
\put(18,-2){\circle*{3}}\put(9,10){\circle*{3}}\put(27,10){\circle*{3}}
\put(0,22){\circle*{3}}\put(18,22){\circle*{3}}       \end{picture}}

\newcommand{\hzero}{\begin{picture}(10,10)(-5,0)\thicklines\put(1.5,2.5){\circle*{2}}\put(-5,0){${\scriptstyle 0}$}\end{picture}}
\newcommand{\hone}{
\begin{picture}(10,10)(-5,0)\thicklines\put(1.5,-2){\line(0,1){10}}
           \put(1.5,-2){\circle*{2}}\put(-5,-3){${\scriptstyle 0}$}\put(1.5,8){\circle*{2}}
           \put(-5,7){${\scriptstyle 1}$}\end{picture}}
\newcommand{\htwotwo}{
\begin{picture}(20,14)(-5,0)\thicklines
          \put(5,-2){\line(-1,2){5}}\put(5,-2){\line(1,2){5}}
           \put(5,-2){\circle*{2}}\put(-2,-4){${\scriptstyle 0}$}
           \put(0,8){\circle*{2}}\put(-6,6){${\scriptstyle 2}$}
           \put(10,8){\circle*{2}}\put(11,6){${\scriptstyle 1}$}
       \end{picture}}
\newcommand{\htwoone}{\begin{picture}(10,20)(-5,0)\thicklines
           \put(1.5,-3){\line(0,1){20}}
           \put(1.5,-2){\circle*{2}}\put(-5,-4){${\scriptstyle 0}$}
           \put(1.5,8){\circle*{2}}\put(-5,6){${\scriptstyle 1}$}
           \put(1.5,18){\circle*{2}}\put(-5,16){${\scriptstyle 2}$}
       \end{picture}}
\newcommand{\hthreesix}{\begin{picture}(22,12)(-1,0)\thicklines
          \put(9,-2){\line(-3,4){9}}\put(9,-2){\line(3,4){9}}\put(9,-2){\line(0,1){12}}
           \put(9,-2){\circle*{2}}\put(2,-4){${\scriptstyle 0}$}
           \put(9,10){\circle*{2}}\put(7,12){${\scriptstyle 2}$}
           \put(0,10){\circle*{2}}\put(-2,12){${\scriptstyle 3}$}
           \put(18,10){\circle*{2}}\put(16,12){${\scriptstyle 1}$}
       \end{picture}}
\newcommand{\hthreethree}{\begin{picture}(26,26)(-3,0)\thicklines
          \put(9,-2){\line(-3,4){9}}\put(9,-2){\line(3,4){9}}\put(18,10){\line(0,1){12}}
           \put(9,-2){\circle*{2}}\put(2,-4){${\scriptstyle 0}$}
           \put(18,22){\circle*{2}}\put(19,22){${\scriptstyle 3}$}
           \put(0,10){\circle*{2}}\put(-6,10){${\scriptstyle 2}$}
           \put(18,10){\circle*{2}}\put(19,10){${\scriptstyle 1}$}
       \end{picture}}
\newcommand{\hthreefour}{\begin{picture}(26,26)(-3,0)\thicklines
          \put(9,-2){\line(-3,4){9}}\put(9,-2){\line(3,4){9}}\put(18,10){\line(0,1){12}}
           \put(9,-2){\circle*{2}}\put(2,-4){${\scriptstyle 0}$}
           \put(18,22){\circle*{2}}\put(19,22){${\scriptstyle 2}$}
           \put(0,10){\circle*{2}}\put(-6,10){${\scriptstyle 3}$}
           \put(18,10){\circle*{2}}\put(19,10){${\scriptstyle 1}$}
       \end{picture}}
\newcommand{\hthreefive}{\begin{picture}(26,26)(-3,0)\thicklines
          \put(9,-2){\line(-3,4){9}}\put(9,-2){\line(3,4){9}}\put(0,10){\line(0,1){12}}
           \put(9,-2){\circle*{2}}\put(2,-4){${\scriptstyle 0}$}
           \put(0,22){\circle*{2}}\put(-6,22){${\scriptstyle 3}$}
           \put(0,10){\circle*{2}}\put(-6,10){${\scriptstyle 2}$}
           \put(18,10){\circle*{2}}\put(19,10){${\scriptstyle 1}$}
       \end{picture}}
\newcommand{\hthreetwo}{\begin{picture}(26,26)(-3,0)\thicklines
          \put(9,10){\line(-3,4){9}}\put(9,10){\line(3,4){9}}\put(9,-2){\line(0,1){12}}
           \put(9,-2){\circle*{2}}\put(2,-4){${\scriptstyle 0}$}
           \put(9,10){\circle*{2}}\put(3,8){${\scriptstyle 1}$}
           \put(0,22){\circle*{2}}\put(-6,22){${\scriptstyle 3}$}
           \put(18,22){\circle*{2}}\put(19,22){${\scriptstyle 2}$}
       \end{picture}}
\newcommand{\hthreeone}{\begin{picture}(12,36)(-7,0)\thicklines
          \put(1.5,-2){\line(0,1){36}}
           \put(1.5,-2){\circle*{2}}\put(-5,-4){${\scriptstyle 0}$}
           \put(1.5,10){\circle*{2}}\put(-5,8){${\scriptstyle 1}$}
           \put(1.5,22){\circle*{2}}\put(-5,20){${\scriptstyle 2}$}
           \put(1.5,34){\circle*{2}}\put(-5,32){${\scriptstyle 3}$}
\end{picture}}
\newcommand{\hfourone}{\begin{picture}(28,28)(-6,0)\thicklines
          \put(9,-2){\line(-3,4){9}}\put(9,-2){\line(3,4){9}}\put(9,-2){\line(0,1){12}}\put(18,10){\line(0,1){12}}
\put(18,22){\circle*{2}}\put(19,24){${\scriptstyle 3}$}
\put(9,-2){\circle*{2}}\put(2,-4){${\scriptstyle 0}$}
\put(9,10){\circle*{2}}\put(7,12){${\scriptstyle 2}$}
\put(0,10){\circle*{2}}\put(-6,12){${\scriptstyle 4}$}
\put(18,10){\circle*{2}}\put(19,12){${\scriptstyle 1}$}
       \end{picture}}

\newcommand{\hfourthree}{\begin{picture}(28,28)(-6,-2)\thicklines
\put(9,-2){\line(-3,4){9}}\put(9,-2){\line(3,4){9}}\put(9,-2){\line(0,1){12}}\put(9,10){\line(0,1){12}}
\put(9,22){\circle*{2}}\put(4,24){${\scriptstyle 3}$}
\put(9,-2){\circle*{2}}\put(2,-4){${\scriptstyle 0}$}
\put(9,10){\circle*{2}}\put(4,12){${\scriptstyle 2}$}
\put(0,10){\circle*{2}}\put(-6,12){${\scriptstyle 4}$}
\put(18,10){\circle*{2}}\put(19,12){${\scriptstyle 1}$}
       \end{picture}}
\newcommand{\hfourfour}{\begin{picture}(26,36)(-3,-5)\thicklines
\put(9,8){\line(-3,4){9}}\put(9,8){\line(3,4){9}}
\put(0,20){\line(0,1){12}}\put(9,-2){\line(0,1){10}}
           \put(9,-2){\circle*{2}}\put(2,-4){${\scriptstyle 0}$}
           \put(9,8){\circle*{2}}\put(2,6){${\scriptstyle 1}$}
           \put(0,32){\circle*{2}}\put(-6,32){${\scriptstyle 4}$}
           \put(0,20){\circle*{2}}\put(-6,20){${\scriptstyle 3}$}
           \put(18,20){\circle*{2}}\put(19,20){${\scriptstyle 2}$}
       \end{picture}}

\newcommand{\tzeroC}[1]{\begin{picture}(5,10)(-1,0)
   {\color{#1}\put(1.5,2.5){\circle*{3}}}
   \end{picture}}
\newcommand{\toneB}{{\color{blue}
\begin{picture}(5,20)(-1,0)\thicklines
           \put(1.5,-2){\line(0,1){20}}
           \put(1.5,-2){\circle*{3}}\put(1.5,18){\circle*{3}}
           \end{picture}}}
\newcommand{\toneR}{{\color{red}
\begin{picture}(5,20)(-1,0)\thicklines
           \put(1.5,-2){\line(0,1){20}}
           \put(1.5,-2){\circle*{3}}\put(1.5,18){\circle*{3}}
           \end{picture}}}

\newcommand{\ttwotwoC}{
\begin{picture}(20,20)(-1,0)\thicklines
      {\color{blue} \put(10,-2){\line(-1,2){10}}\put(0,18){\circle*{3}}}{\color{green}
\put(10,-2){\line(1,2){10}}\put(20,18){\circle*{3}}}     
\put(10,-2){\circle*{3}}
       \end{picture}}
\newcommand{\ttwotwoBR}{
\begin{picture}(20,20)(-1,0)\thicklines
 {\color{blue}
 \put(10,-2){\line(-1,2){10}}\put(0,18){\circle*{3}}
 }{\color{red}
 \put(10,-2){\line(1,2){10}}\put(20,18){\circle*{3}}
 }\put(10,-2){\circle*{3}}
       \end{picture}}

\newcommand{\ttwooneG}{{\color{green}
\begin{picture}(5,20)(-1,0)\thicklines
           \put(1.5,-3){\line(0,1){20}}
           \put(1.5,-2){\circle*{3}}\put(1.5,8){\circle*{3}}\put(1.5,18){\circle*{3}}
       \end{picture}}}
\newcommand{\tthreefiveC}{
 \begin{picture}(20,20)(-1,0)\thicklines
  {\color{blue} \put(10,-2){\line(-1,2){10}}\put(0,18){\circle*{3}}}{\color{green}
\put(10,-2){\line(0,1){20}}\put(10,18){\circle*{3}}}{\color{red}  
\put(10,-2){\line(1,2){10}}\put(20,18){\circle*{3}}}{\color{black} \put(10,-2){\circle*{3}}}
 \end{picture}}
\newcommand{\tthreethreeBRG}{
 \begin{picture}(20,24)(-2,0)\thicklines
{\color{blue} \put(9,-2){\line(-3,4){9}}\put(0,10){\circle*{3}}}{\color{green}
\put(18,10){\line(0,1){12}}\put(18,22){\circle*{3}}}{\color{red}  
\put(9,-2){\line(3,4){9}}}{\color{black} \put(9,-2){\circle*{3}} \put(18,10){\circle*{3}}}
 \end{picture}}
\newcommand{\tthreethreeGRB}{
 \begin{picture}(20,24)(-2,0)\thicklines
 {\color{green}  \put(9,-2){\line(-3,4){9}}\put(0,10){\circle*{3}}}{\color{blue}
\put(18,10){\line(0,1){12}}\put(18,22){\circle*{3}}}{\color{red}  
\put(9,-2){\line(3,4){9}}}{\color{black} \put(9,-2){\circle*{3}} \put(18,10){\circle*{3}}}
 \end{picture}} 
\newcommand{\tthreethreeRG}{
 \begin{picture}(23,24)(-2,0) \thicklines
 {\color{red}
 \put(9,-2){\line(-3,4){9}}\put(0,10){\circle*{3}}
 }{\color{green}
 \put(18,10){\line(0,1){12}}\put(18,22){\circle*{3}}\put(9,-2){\line(3,4){9}}\put(18,10){\circle*{3}}
   } \put(9,-2){\circle*{3}}
 \end{picture}}
 \newcommand{\tthreethreeBG}{
 \begin{picture}(23,24)(-1,0) \thicklines
 {\color{blue}
 \put(9,-2){\line(-3,4){9}}\put(0,10){\circle*{3}}
 }{\color{green}
 \put(18,10){\line(0,1){12}}\put(18,22){\circle*{3}}\put(9,-2){\line(3,4){9}}\put(18,10){\circle*{3}}
   } \put(9,-2){\circle*{3}}
 \end{picture}}
\newcommand{\tthreetwoC}{
  \begin{picture}(20,24)(-1,0) \thicklines
{\color{blue}\put(9,10){\line(-3,4){9}}\put(0,22){\circle*{3}}}{\color{green}\put(9,10){\line(3,4){9}}\put(18,22){\circle*{3}}}{\color{red}\put(9,-2){\line(0,1){12}}\put(9,-2){\circle*{3}}}{\color{black}\put(9,10){\circle*{3}}}
   \end{picture}}

\newcommand{\tfouroneC}{
 \begin{picture}(23,24)(-2,0) \thicklines
  {\color{blue}
  \put(9,-2){\line(-3,4){9}}\put(0,10){\circle*{3}}
  }{\color{red}
  \put(9,-2){\line(0,1){12}}\put(9,10){\circle*{3}}
  }{\color{green}
  \put(9,-2){\line(3,4){9}}\put(18,10){\line(0,1){12}}\put(18,22){\circle*{3}}\put(18,10){\circle*{3}}
   }\put(9,-2){\circle*{3}}
 \end{picture}}

\begin{document}

\title[Hopf algebras of permutations and trees]{Cocommutative Hopf algebras \\of
       permutations and trees}
\author{Marcelo Aguiar}
\address{Department of Mathematics\\ Texas A\&M University\\
 College Station, TX 77843, USA}
 \email{maguiar@math.tamu.edu}
 \urladdr[Marcelo Aguiar]{http://www.math.tamu.edu/$\sim$maguiar}

\author{Frank Sottile}
\address{Department of Mathematics\\
         Texas A\&M University\\
         College Station, TX 77843, USA}
\email{sottile@math.tamu.edu}
\urladdr{http://www.math.tamu.edu/$\sim${}sottile}

\thanks{Aguiar supported in part by NSF grant DMS-0302423}
\thanks{Sottile supported in part by NSF CAREER grant DMS-0134860, the
  Clay Mathematics Institute, and MSRI}
  \thanks{We thank Lo\"\i c Foissy and Ralf Holtkamp for interesting comments and remarks.}

\begin{abstract} Consider the coradical filtrations of the Hopf algebras of planar binary
trees of Loday and Ronco and of permutations of Malvenuto and Reutenauer.  We
give explicit isomorphisms showing that the associated graded Hopf algebras are
dual to the cocommutative Hopf algebras introduced in the late 1980's by
Grossman and Larson. These Hopf algebras are constructed from ordered trees and
heap-ordered trees, respectively. These results follow from the fact that
whenever one starts from a Hopf algebra that is a cofree  graded  coalgebra, the
associated graded Hopf algebra is a shuffle Hopf algebra. 
\end{abstract}

\maketitle

\section*{Introduction}\label{S:intro}
In the late 1980's, Grossman and Larson constructed  
several
cocommutative
Hopf algebras from 
different
families of trees (rooted, ordered,
heap-ordered), in connection to the symbolic algebra of differential operators~\cite{GL89,GL90}.
Other Hopf algebras of trees have arisen lately in a variety of
contexts, including
 the Connes-Kreimer Hopf algebra in renormalization theory~\cite{CK} and the Loday-Ronco
Hopf algebra in the theory of associativity breaking~\cite{LR98,LR02}. The latter is
closely related to other important Hopf algebras in algebraic combinatorics, including
the Malvenuto-Reutenauer Hopf algebra~\cite{MR95} and the Hopf algebra
of quasi-symmetric functions~\cite{Mal, Re93, St99}. 

This universe of Hopf algebras of trees is summarized below.

\smallskip
\begin{center}
\begin{tabular}{c|c|c}
 & \bf{Family of trees} & \bf{Hopf algebra}  \\
\hline\hline
  & rooted trees &  \\
\cline{2-2}
 Grossman-Larson & ordered trees   & non-commutative, \\
\cline{2-2}
89-90 & heap-ordered  & cocommutative \\
 & trees &  \\
  \hline\hline
 Loday-Ronco  & planar binary  &  non-commutative,\\
 98 & trees&  non-cocommutative\\
\hline\hline
Connes-Kreimer & rooted trees &  commutative, \\
  98    & &  non-cocommutative\\
 \hline\hline
\end{tabular}
\end{center}
\medskip

Recent independent work of Foissy~\cite{Foi02,Foi02b} and 
Hoffman~\cite{Ho02} showed that the Hopf algebra of
Connes-Kreimer is dual to the Hopf algebra of rooted trees of Grossman-Larson. 
This Hopf algebra also arises as the universal enveloping algebra
of the free {\em pre-Lie algebra} on one generator, viewed as a
Lie algebra~\cite{CL02}.
Foissy~\cite{Foi02b} and Holtkamp~\cite{Ho03} showed that the Hopf algebra of Connes-Kreimer is a 
quotient of the Hopf algebra of Loday-Ronco, 
see also~\cite{ASb}.

We give explicit isomorphisms which show that the Grossman-Larson Hopf algebras of ordered trees and of
heap-ordered trees are dual to the {\em associated graded}\/ Hopf algebras to the Hopf algebra
$\YSym$ of planar binary trees of Loday and Ronco and the Hopf algebra $\SSym$ of
permutations of Malvenuto and Reutenauer, respectively. This is done in
Theorems~\ref{T:GL-YSym} and~\ref{T:HGL-SSym}. The case of heap-ordered trees requires the assumption that the base field be of characteristic $0$.
 We establish this case in Section~\ref{S:SSym} by making use of the {\em first Eulerian idempotent}.

The essential tool we use is the {\em monomial} basis of $\YSym$ and $\SSym$ introduced in
our previous works~\cite{ASa,ASb}. The
explicit isomorphisms are in terms of the dual bases of ordered and heap-ordered trees of
Grossman-Larson and of the  monomial bases of $\YSym$ and $\SSym$, respectively. These results
provide unexpected combinatorial descriptions for the associated graded Hopf algebras to
$\YSym$ and $\SSym$. On the other hand, together with  the result of Foissy and Hoffman,
they connect all Grossman-Larson Hopf algebras of trees to the mainstream of combinatorial
Hopf algebras.

It follows from our results that the associated graded Hopf algebras to the Hopf
algebras of Loday-Ronco and Malvenuto-Reutenauer are commutative, a fact which is not
obvious from the explicit description of the product of these algebras.  
Greg Warrington noticed this for the Malvenuto-Reutenauer Hopf algebra  and Lo\"\i c Foissy made us aware that
the associated graded Hopf algebra to any graded connected Hopf algebra is
always commutative (private communications). A related well-known
fact is that the associated graded Hopf algebra to a {\em cofree} graded connected Hopf algebra is a shuffle Hopf algebra. We recall these and related results with their proofs in Section~\ref{S:cofree}.
This also implies that the algebras of Grossman and Larson
are  tensor Hopf algebras (Corollaries~\ref{C:free-ordered}
and~\ref{C:free-heap}). It was known from~\cite{GL89} that these
algebras are free.


\section{Cofree  graded coalgebras and Hopf algebras}\label{S:cofree}

A coalgebra $(C,\Delta,\epsilon)$ over a field $\field$ is called {\em graded} if there is
given a decomposition $C=\oplus_{k\geq 0}C^k$ of $C$ as a direct sum of $\field$-subspaces
$C^k$ such that
\[\Delta(C^k)\subseteq\sum_{i+j=k}C^i\otimes C^j \text{ \ and  \ }
\epsilon(C^k)=0\ \ \forall\, k\neq 0\,.\]
The coalgebra is said to be {\em graded connected} if in addition $C^0\cong\field$.

\begin{defi}\label{D:cofree} A graded coalgebra $Q=\oplus_{k\geq 0}Q^k$
is said to be {\em cofree} if it satisfies the following universal property.
Given a graded coalgebra $C=\oplus_{k\geq 0}C^k$ and a linear map
$\varphi:C\to Q^1$ with $\varphi(C^k)=0$ when $k\neq 1$,
there is a unique morphism of graded coalgebras
$\hat{\varphi}:C\to Q$ such that the following diagram commutes
\[\xymatrix{{\ C\ }\ar@{-->}[rr]^{\hat{\varphi}}\ar[dr]_{\varphi} &
   &{Q}\ar[ld]^{\pi}\\ & {Q^1} }\]
where $\pi:Q\to Q^1$ is the canonical projection.
\end{defi}

Let $V$ be a vector space and set
 \[Q(V)\ :=\ \bigoplus_{k\geq 0} V^{\ten k}\,.\]
We write 
elementary
tensors from $V^{\ten k}$ as
$x_1\iten x_2\iten\cdots\iten x_k$ ($x_i\in V$) and identify $V^{\ten 0}$ with
 $\field$. The space $Q(V)$, graded by $k$, becomes a graded connected coalgebra with the
\emph{deconcatenation} coproduct
\begin{equation}\label{E:deconcat}
\Delta(x_1\iten x_2\iten\cdots\iten x_k)\ =\
   \sum_{i=0}^k\
       (x_1\iten\cdots\iten x_i)\ten(x_{i+1}\iten\cdots\iten x_k)
\end{equation}
and counit given by projection onto $V^{\ten 0}=\field$. Moreover, $Q(V)$ is a cofree graded coalgebra~\cite[Lemma 12.2.7]{Swe}. It is in fact graded connected.

By universality, any cofree graded coalgebra $Q$ is isomorphic to $Q(V)$, where
$V=Q^1$. We refer to $Q(V)$ as the cofree graded coalgebra
{\em cogenerated} by $V$.

\begin{rem} The functor $Q$ from vector spaces to graded coalgebras is right
adjoint to the forgetful functor $C\mapsto C^1$ from graded coalgebras to vector
spaces. $Q(V)$ is {\em not} cofree in the  category of all coalgebras over $\field$.
However, $Q(V)$ is still cofree in the  category of connected coalgebras
in the sense of Quillen~\cite[Appendix B, Proposition 4.1]{Qui}. See also~\cite[Theorem 12.0.2]{Swe}.

\end{rem}

We are interested in Hopf algebra structures on cofree graded
coalgebras. There is recent important work of Loday and Ronco in
this direction~\cite{LR04}, but their results are not prerequisites
for our work.

In the classical Hopf algebra literature usually only one Hopf
algebra structure on $Q(V)$ is considered: the shuffle Hopf algebra.
It is well-known that this is the {\em only} Hopf algebra structure on $Q(V)$ for which the algebra structure
 preserves the grading; this may  be
deduced from~\cite[Theorem 12.1.4]{Swe} but we provide a direct proof below (Proposition~\ref{P:shuffle}).  
There are, however, many  naturally occurring Hopf algebras that are cofree
graded coalgebras and for which the algebra structure {\em does not}
 preserve the grading; see Examples~\ref{Ex:cofree}.

{\em The shuffle Hopf algebra}. Let $V$ be an arbitrary vector space. There is
an algebra structure on $Q(V)$ defined recursively by 
\[x\cdot 1=x=1\cdot x\]
for $x\in V$, and
\begin{multline*}
(x_1\iten\cdots\iten x_j)\cdot(y_1\iten\cdots\iten y_k)=\\
x_1\iten\Bigl((x_2\iten\cdots\iten x_j)\cdot(y_1\iten\cdots\iten y_k)\Bigr)+
y_1\iten\Bigl((x_1\iten\cdots\iten x_j)\cdot(y_2\iten\cdots\iten y_k)\Bigr)\,.
\end{multline*}
Together with the graded coalgebra structure~\eqref{E:deconcat}, this gives a Hopf
algebra which is denoted $\Sh(V)$ and called the {\it shuffle Hopf algebra of $V$}.

A Hopf algebra $H$ is called graded if it is a graded coalgebra and the
multiplication and unit preserve the grading: 
\[H^j\cdot H^k\subseteq H^{j+k}\,,\ \ 1\in H^0\,.\]
The shuffle Hopf algebra $\Sh(V)$ is a graded Hopf algebra. As mentioned, it is the only
such structure that a cofree graded coalgebra admits.

\begin{prop} \label{P:shuffle}
Let  $H=\oplus_{k\geq 0}H^k$ be a graded Hopf algebra which is cofree as a graded coalgebra.
Then there is an isomorphism of graded Hopf algebras 
 \[H\cong\Sh(H^1)\,.\]
\end{prop}

\begin{proof}  We may assume that $H=Q(V)$, with $V=H^1$. By hypothesis, the
multiplication map is a morphism of graded Hopf algebras $m:H\ten H\to H$, where the
component of degree $k$ of $H\ten H$ is
$\sum_{i+j=k}H^i\ten H^j$. By cofreeness, $m$ is uniquely determined by
the  composite 
\[H\ten H\map{m}H\map{\pi}H^1\,,\]
which in turn reduces to
\[(H\ten H)^1=H^0\ten H^1+H^1\ten H^0\map{m}H^1\,.\]
Also by hypothesis, $H^0=\field\cdot 1$ where $1$ is the unit element of $H$.
Hence the above map, and then $m$, are determined by
\[1\ten x\mapsto x \text{ \ and \ }x\ten 1\mapsto x\,.\]
This shows that there is a unique multiplication on $H$ that makes it a graded Hopf
algebra. 
Since the multiplication of the shuffle Hopf algebra of $H^1$ is one such map,
it is the only one. Thus, $H$ is the shuffle Hopf algebra of $H^1$.
\end{proof}

\medskip

{\em The tensor Hopf algebra}. Let $V$ be a vector space and set
 \[T(V)\ :=\ \bigoplus_{k\geq 0} V^{\ten k}\,.\]
As a vector space, $T(V)=Q(V)$. The space $T(V)$ becomes a graded algebra
under the {\em concatenation product}
\[ (x_1\iten\cdots\iten x_i)\cdot(y_{1}\iten\cdots\iten y_j)=
x_1\iten \cdots \iten x_i\iten y_1\iten\cdots\iten y_j\]
and unit $1\in V^{\ten 0}=\field$. Moreover, $T(V)$ is the free algebra on $V$.

If $V$ is finite dimensional, the graded dual of $Q(V)$ is the tensor algebra  $T(V^*)$.

There is a graded Hopf algebra structure on $T(V)$ uniquely determined by
\[\Delta(x)=1\ten x+ x\ten 1 \text{ \and \ } \epsilon(x)=0\]
for $x\in V$. This the tensor Hopf algebra. An argument dual to that of
Proposition~\ref{P:shuffle} shows that it is the only
graded Hopf algebra structure that a free algebra admits.

\medskip

{\em The coradical filtration}. Let $H=\oplus_{k\geq 0}H^k$ be a  Hopf algebra that is graded as a coalgebra. We do not insist that the algebra structure of $H$ preserves this grading. Let $F^0(H):=H^0$ and let $ F^k(H)$ consist of those elements $h\in H$ such that in the iterated coproduct $\Delta^{(k)}(h)$ every term has a tensor factor from $F^0(H)$. It follows that $F^k(H)\subseteq F^{k+1}(H)$ and
$H^k\subseteq F^k(H)$.

Suppose $H$ is connected, i.e., $F^0(H)=H_0=\field$. In this case,  $F^0(H)$ is the {\em coradical} of $H$ and the subspaces $F^k(H)$ form the {\em coradical filtration} of $H$~\cite[Chapter 5]{Mo93}. It is known that  
\[
H=\bigcup_{k\geq 0}F^k(H)\,, \quad
\Delta\bigl(F^k(H)\bigr)\subseteq \sum_{i+j=k}F^i(H)\otimes F^j(H)\,,
\text{ \ and \ }
F^j(H)\cdot F^k(H)\subseteq F^{j+k}(H)\,.
\]
These results hold in greater generality; see~\cite[Theorem 5.2.2, Lemma 5.2.8]{Mo93}.

Let $\gr(H)$ be the  graded Hopf algebra associated to the coradical filtration, 
\[\gr(H)=F^0(H)\oplus F^{1}(H)/F^0(H)\oplus\cdots\oplus F^{k+1}(H)/F^k(H)\oplus\cdots\]
If $m$ and $\Delta$ are the operations of $H$, then the
operations of $\gr(H)$ are induced by the compositions
\begin{gather*}
\xymatrix{ {F^j(H)\otimes F^k(H)}\ar[r]^-{m} &
              {F^{j+k}(H)}\ar@{->>}[r] & {F^{j+k}(H)/F^{j+k-1}(H)} }\,,\\
 \xymatrix{  
              {F^{k}(H)}\ar[r]^-{\Delta} & {\sum_{i+j=k}F^i(H)\otimes F^j(H)}\ar@{->>}[r] & {\sum_{i+j=k} F^{i+j}(H)/F^{i+j-1}(H)} } \,.
\end{gather*}
The main goal of this paper is to obtain explicit combinatorial descriptions for the
associated graded Hopf algebras to the Hopf algebras $\YSym$ and $\SSym$ of
Examples~\ref{Ex:cofree}. This is done in Sections~\ref{S:YSym} and~\ref{S:SSym}. These Hopf algebras are cofree graded coalgebras, so we discuss the coradical filtration for
such Hopf algebras first.

Let $H=Q(V)$ be a Hopf algebra that is a cofree graded  coalgebra.
We have $H^0=\field$, $H^1=V=P(H)$, the space of primitive elements
of $H$, and $H^k=V^{\ten k}$. As before,
we do not require that the algebra structure of $H$ preserves this grading.
It is easy to see that
\[F^k (H)\ :=\ H^0 \oplus H^1 \oplus \cdots \oplus H^k\,.\]
Therefore, $\gr(H)\cong H$ as graded coalgebras canonically, and the multiplication has been altered
by removing terms of lower degree from a homogeneous product.
More precisely, 
if $m$ is the multiplication map on $H$, then the
multiplication on $\gr(H)$ is  the composition
\[
\xymatrix{ {H^j\otimes H^k}\ar[r]^-{m} &
              {F^{j+k}(H)}\ar@{->>}[r] & {H^{j+k}} }\,.\]

\begin{prop}\label{P:shuffle-2}
Let $H$ be a Hopf algebra that is a cofree graded  coalgebra. 
Then its associated graded Hopf algebra $\gr(H)$ is the shuffle Hopf algebra $\Sh(H^1)$. 
In particular, $\gr(H)$  is commutative.
\end{prop}
\begin{proof} Since $H\cong\gr(H)$ as graded coalgebras, Proposition~\ref{P:shuffle} applies to $\gr(H)$.
\end{proof}

The commutativity of the associated graded Hopf algebra holds in greater
generality. The following result was pointed out to us by Foissy.

\begin{prop} \label{P:graded-comm}
Let  $H$ be a graded connected Hopf algebra. Then $\gr(H)$ is commutative.
\end{prop}
\begin{proof} We show that  $[F^j(H),F^k(H)]\subseteq F^{j+k-1}(H)$, and hence commutators vanish in $\gr(H)$.

 It follows from the definition of the coradical filtration that for any $h\in F^a(H)$ every term in $\Delta^{(a+b-1)}(h)$ contains at least $b$ factors from $F^0(H)=\field$.

Let $x\in F^j(H)$ and $y\in F^k(H)$. Every term in $\Delta^{(j+k-1)}(x)$ contains at least $k$ factors from $\field$ and every term in $\Delta^{(j+k-1)}(y)$ contains at least $j$ factors from $\field$. Write
\[\Delta^{(j+k-1)}(x)=\sum x_1\otimes\cdots\otimes x_{j+k} \text{ \ and \ }
\Delta^{(j+k-1)}(y)=\sum y_1\otimes\cdots\otimes y_{j+k}\,.\]
 Consider those terms in 
\[\Delta^{(j+k-1)}(xy)=\Delta^{(j+k-1)}(x)\Delta^{(j+k-1)}(y)=
\sum x_1y_1\otimes\cdots\otimes x_{j+k}y_{j+k}
\]
in which {\em none} of the $j+k$ factors are from $\field$. By the pigeon-hole principle, these terms must be such that for each $i=1,\ldots, j+k$ either $x_i\in\field$ or $y_i\in\field$. Therefore, these terms satisfy
\[x_1y_1\otimes\cdots\otimes x_{j+k}y_{j+k}=y_1x_1\otimes\cdots\otimes y_{j+k}x_{j+k}\,.\]
The right-hand side is a term in $\Delta^{(j+k-1)}(yx)$, and by symmetry this gives all
terms in $\Delta^{(j+k-1)}(yx)$ in which none of the factors are from $\field$.
These cancel in $\Delta^{(j+k-1)}(xy-yx)$. Thus,
every term in $\Delta^{(j+k-1)}(xy-yx)$ contains at least one factor from $\field$, which proves that $xy-yx\in F^{j+k-1}(H)$.
\end{proof}

\begin{rem} Consider the coradical filtration of an arbitrary (not necessarily
  graded or connected) Hopf algebra. The same argument as above shows that if
  the coradical $F^0(H)$ lies in the center of $H$,  then the associated graded Hopf
  algebra is commutative. 
\end{rem}

\medskip

The cofree graded coalgebras we are interested in carry a second grading. With respect
to this second grading, but not with respect to the original one, they  are in fact graded
Hopf algebras. The general setup is as follows.

Suppose $V=\oplus_{i\geq 1}V_i$ is a graded space and each $V_i$ is finite dimensional. 
Then $Q(V)$ carries another grading, for which the elements of
$V_{i_1}\ten\cdots\ten V_{i_k}$ have degree $i_1+\cdots+i_k$. In this situation, we 
refer to $k$ as the {\em length} and to $i_1+\cdots+i_k$ as the {\em weight}. The
homogeneous components of the two gradings on $Q(V)$ are thus
\[Q(V)^k:=V^{\ten k} \text{ \ and \ }
Q(V)_n:=\bigoplus_{\substack{k\geq 0\\i_1+\cdots+i_k=n}}V_{i_1}\ten\cdots\ten V_{i_k}\,.\]
Note that each $Q(V)_n$ is finite dimensional.
Let $V^*:=\oplus_{i\geq 1}V_i^*$  denote the {\em graded dual} of $V$.
The graded dual of $Q(V)$ with
respect to the grading by weight is the tensor algebra
$T(V^*)$, and the graded dual of $\Sh(V)$ with respect to the grading
by weight is the tensor Hopf algebra $T(V^*)$.

\begin{exas}\label{Ex:cofree}
 We give some examples of cofree graded coalgebras.

(1) {\em The Hopf algebra of  quasi-symmetric functions.} This Hopf algebra, often
denoted $\QSym$, has a linear basis $M_\alpha$ indexed by compositions 
$\alpha=(a_1,\ldots,a_k)$
(sequences of positive integers). See~\cite{Mal, Re93, St99} for more details.
$\QSym$ is a cofree graded coalgebra, as follows.  Let $V$ be the subspace
linearly spanned by the elements
$M_{(n)}$, $n\geq 1$.  Then 
$\QSym\cong Q(V)$ via 
\[M_{(a_1,\ldots,a_k)}\longleftrightarrow M_{(a_1)}\iten\cdots\iten M_{(a_k)}\,.\]
This isomorphism identifies $V^{\ten k}$ with the subspace of $\QSym$ spanned by the
elements $M_\alpha$ indexed by compositions of length $k$.  
$\QSym$ is not a shuffle Hopf algebra: the product does not preserve the
grading by length. For instance,
\[M_{(n)}\cdot M_{(m)}=M_{(n,m)}+M_{(m,n)}+M_{(n+m)}\,.\]
In this case, $V$ is graded by $n$, and the grading by weight assigns degree
$a_1+\dotsb+a_k$ to $M_\alpha$. The Hopf algebra structure of $\QSym$ does preserve
the grading by weight. 

This is an example of a {\em quasi-shuffle} Hopf algebra~\cite{Ho00,Ha}.
According to~\cite[Theorem 3.3]{Ho00}, any (commutative) quasi-shuffle 
Hopf algebra is  isomorphic to a shuffle Hopf algebra. 
The isomorphism does not however preserve the grading by length, 
and thus its structure as a cofree graded coalgebra.
For more on the cofreeness of $\QSym$, see~\cite[Theorem 4.1]{ABS}.

(2) {\em The Hopf algebra of planar binary trees.} This Hopf algebra was introduced by Loday
and Ronco~\cite{LR98,LR02}. We denote it by $\YSym$. It is known that $\YSym$ is a cofree
graded coalgebra~\cite[Theorem 7.1, Corollary 7.2]{ASb}. The product of $\YSym$
does not preserve the grading by length (but it preserves the grading by weight). 
$\YSym$ is not a shuffle Hopf algebra,
 not even a quasi-shuffle Hopf algebra. See Section~\ref{S:YSym} for more details.

(3) {\em The Hopf algebra of permutations.} This Hopf algebra was introduced by Malvenuto
and Reutenauer~\cite{Mal,MR95}. We denote it by $\SSym$. As for $\YSym$, $\SSym$ is
a cofree  graded coalgebra~\cite[Theorem 6.1, Corollary 6.3]{ASa} and is neither a shuffle
nor quasi-shuffle Hopf algebra. See Section~\ref{S:SSym} for more details.

(4) {\em The Hopf algebra of peaks.} This Hopf algebra was introduced by
Stembridge~\cite{Ste} and is often denoted $\Pi$. It has a linear basis indexed by {\em odd
compositions} (sequences of non-negative odd integers). It has been recently shown that
$\Pi$ is a cofree graded coalgebra~\cite[Theorem 4.3]{Hsi},
see also~\cite[Proposition 3.3]{Sch}.
\end{exas}

\section{The Hopf algebra of ordered trees} \label{S:YSym}

We show that the graded dual to $\gr(\YSym)$ is isomorphic to the cocommutative
Hopf algebra of ordered trees defined by Grossman and Larson~\cite{GL89}.

We first review the definition of the Hopf algebra of ordered trees.

For the definition of ordered trees (also called rooted planar trees), see~\cite[page 294]{St86}.
The ordered trees with 1, 2, 3, and 4 nodes
are shown below:
\[
  \tzero,\quad\qquad\tone,\quad\qquad\ttwoone,\ttwotwo,\quad\qquad
  \tthreeone,\ \tthreetwo,\ \tthreethree,\ \tthreefour,\ \tthreefive\,.
\]
Given two ordered trees $x$ and $y$, we may join them together at their
roots to obtain another ordered tree $x\iten y$, where the nodes
 of $x$ are to the left of those of $y$:
\[\epsfxsize=0.4in\epsfbox{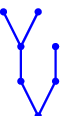}\ \epsfxsize=0.35in\epsfbox{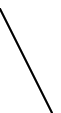}\
 \epsfxsize=0.3in\epsfbox{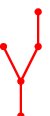}\ \
 \raisebox{24pt}{=}\ \epsfxsize=0.9in\epsfbox{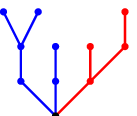}\ .\]

An ordered tree is {\em planted} if its root has a unique child.
Every ordered tree $x$ has a unique decomposition
\begin{equation}\label{E:planted}
x=x_1\iten\cdots\iten x_k
\end{equation}
 into planted trees $x_1,\dotsc,x_k$, corresponding to the branches at the root of $x$.
These are the {\em planted components} of $x$.

The set of nodes of an ordered tree $x$ is denoted by $\Nod(x)$.
Let $x$ be an ordered tree and $x_1,\dotsc,x_k$ its  planted components, listed from left to
right and (possibly) with multiplicities.
Given a function $f:[k]\to\Nod(y)$ from the set $[k]=\{1,\ldots,k\}$  to the set of nodes of
another ordered tree
$y$, form a new ordered tree $x\tof y$ by
identifying the
root of each component $x_i$ of $x$ with the corresponding node $f(i)$ of $y$.
For this to be an ordered tree, retain the order of any components of $x$ attached to
the same node of $y$, and place them to the left of any children of that node in
$y$. Given a subset $S\subseteq [k]$, say $S=\{i_1<\cdots<i_p\}$, let
\[x_S:=x_{i_1}\iten\cdots\iten x_{i_p}\,.\]
Equivalently, $x_S$ is the tree obtained  by erasing the branches at the root of $x$ which
are not indexed by $S$. Let $S^c=[k]\setminus S$.

\begin{defi}\label{D:GL}
 The  Grossman-Larson Hopf algebra $\GL$ of ordered trees is the
 formal linear span of all ordered trees with product and coproduct as follows.
 Given ordered trees $x$ and $y$ as above, we set
 \begin{eqnarray*}
   x\cdot y &=& \sum_{f:[k]\to\Nod(y)}\!\!\!\! x\tof y\,,\\
  \Delta(x) &=& \sum_{S\subseteq[k]}x_S\otimes x_{S^c}\,,\rule{0pt}{15pt}
 \end{eqnarray*}
the first sum is over all functions  from $[k]$ to the set of nodes of $y$ and
the second is over all  subsets of $[k]$.
$\GL$ is a graded Hopf algebra, where the degree of an ordered tree is one less than
the number of nodes~\cite[Theorem 3.2]{GL89}.
\end{defi}

We give some examples, using colors 
to indicate how the operations are performed 
(they are not part of the structure of an ordered tree).

\[
  \ttwotwoC\cdot\toneR\ \  =\  \tthreetwoC+\tthreethreeBRG+\tthreethreeGRB
      +\tthreefiveC\ =\ \tthreetwo+2\cdot\tthreethree+\tthreefive\,.
\]
\begin{eqnarray*}
  \Delta\bigl(\tfouroneC\,\bigr)\ &=&\
   \tzero\ten\tfouroneC\ +\ \toneB\ten\tthreethreeRG\ +\ 
    \toneR\ten\tthreethreeBG\ +\ \ttwooneG\ten\ttwotwoBR \\
   &&+\ \tfouroneC\ten\tzero\ +\ \tthreethreeRG\ten\toneB\ +\
    \tthreethreeBG\ten\toneR\ +\ \ttwotwoBR\ten\ttwooneG \\
     &=&\ \tzeroC{black}\ten\tfourone\ +\
2\cdot\tone\ten\tthreethree\ +\ \ttwoone\ten\ttwotwo\\
&& +\ \tfourone\ten\tzeroC{black}\ +\ 2\cdot\tthreethree\ten\tone
+\ \ttwotwo\ten\ttwoone.
\end{eqnarray*}

The definition implies 
that $\GL$ is cocommutative and that each planted tree is a
primitive element in $\GL$. (There are 
other primitive elements. In fact, $\GL$ is
isomorphic to the tensor Hopf algebra on the subspace spanned by the set of planted trees. See
Corollary~\ref{C:free-ordered}.)

\smallskip
We follow the notation and terminology of~\cite{ASb}  
for planar binary trees and
the Loday-Ronco Hopf algebra $\YSym$ (much of which is
based on the constructions of~\cite{LR98,LR02}). 

Ordered trees  are in bijection with planar binary trees. Given a planar binary tree $t$,
draw a node on each of its leaves, then collapse all edges of the form $\slash$. The
resulting planar graph, rooted at the node coming from the right-most leaf of $t$, is an ordered tree.
This defines a  bijection $\psi$ from planar
binary trees with $n$ leaves to ordered trees with $n$ nodes.

We will make use of a recursive definition of $\psi$.
Recall the operation
$s\iten t$ between planar binary trees, which is obtained by identifying the
right-most leaf of $s$ with the root of $t$ (putting $s$ under $t$). For instance,
\begin{center}
 \epsfysize=40pt\epsffile{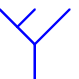}
 \raisebox{12pt}{\Huge $\backslash$}
 \epsfysize=40pt\epsffile{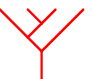}
  \raisebox{18pt}{\large \ =\ \ }
 \epsfysize=40pt\epsffile{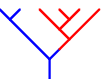}
\end{center}

This operation is associative and so 
any planar binary tree $t$ has a unique maximal decomposition
\begin{equation}\label{E:spikes}
t=t_1\iten t_2\iten\cdots\iten t_k
\end{equation}
in which each $t_i$ is $\iten$-irreducible. Note that a planar binary tree $t$ is $\iten$-irreducible
precisely when it is of the form
\begin{equation}\label{E:nospikes}
t \ = \ \raisebox{-25pt}{\begin{picture}(80,60) \thicklines
 \put(40,0){\line(0,1){20}}\put(40,20){\line(-1,1){40}}\put(40,20){\line(1,1){40}}
 \put(33,28){\line(1,1){32}}
 \put(13,47){$\cdots t'\cdots$}
\end{picture}}
\end{equation}
for some  planar binary tree $t'$ with one less leaf than $t$.

The bijection $\psi$ may be computed
 recursively as follows. First, for $t$ as in~\eqref{E:spikes},
\[\psi(t)=\psi(t_1)\iten \psi(t_2)\iten\cdots\iten \psi(t_k)\,.\]
Second, for $t$ as in~\eqref{E:nospikes}, $\psi(t)$ is obtained by adding a new root
to the ordered tree $\psi(t')$:
\[t \ = \ \raisebox{-25pt}{\begin{picture}(80,60) \thicklines
 \put(40,0){\line(0,1){20}}\put(40,20){\line(-1,1){40}}\put(40,20){\line(1,1){40}}
 \put(33,28){\line(1,1){32}}
 \put(13,48){$\cdots t'\cdots$}
\end{picture}}\ \Rightarrow\ 
\psi(t)\ =\ \raisebox{-25pt}{\begin{picture}(80,60)
\thicklines
 \put(40,0){\line(0,1){20}}\put(40,20){\line(-1,1){40}}\put(40,20){\line(1,1){40}}
 \put(12,48){$\cdots \psi(t')\cdots$}\put(0,60){\circle*{3}}\put(80,60){\circle*{3}}
 {\color{red} \put(40,0){\line(0,1){20}}\put(40,0){\circle*{3}}\put(40,20){\circle*{3}}
 }
\end{picture}}\]
 Finally, $\psi(|)=\tzero\,$ is the unique  ordered tree with one
node. For instance,
\begin{align*}
\psi(\raisebox{-2pt}{\epsffile{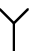}}) &=\tone,\\
  \psi(\epsffile{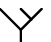})&=
    \psi(\raisebox{-2pt}{\epsffile{figures/1.eps}})\backslash\psi(\raisebox{-2pt}{\epsffile{figures/1.eps}})
    \ =\  \ttwotwo,\\
   \psi(\epsffile{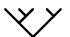})&= \psi(\raisebox{-2pt}{\epsffile{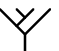}})
     \backslash \psi(\raisebox{-2pt}{\epsffile{figures/1.eps}})\ =\ 
\tthreetwo \ \raisebox{4pt}{\Large$\backslash$}\ \tonebig =\tfourthree\,.
     \end{align*}
Note that $\psi$ identifies $\iten$-irreducible planar binary trees with planted ordered trees.

\smallskip
In~\cite{ASb}, we introduced a linear basis
$M_t$ of $\YSym$, indexed by planar binary trees $t$, which is obtained from the original
basis of Loday and Ronco by a process of M\"obius inversion. We showed that $\YSym$ is a 
cofree graded  coalgebra and the space $V$ of primitive elements is
the linear span of the elements $M_t$ for $t$ a $\iten$-irreducible planar binary
trees~\cite[Theorem 7.1, Corollary 7.2]{ASb}. The isomorphism $Q(V)\cong\YSym$ is
\[M_{t_1}\iten\cdots\iten M_{t_k}\longleftrightarrow M_{t_1\iten\cdots\iten t_k}\,.\]
The resulting grading by length on $\YSym$ is given by the number of $\iten$-irreducible components
in the decomposition of a planar binary tree $t$ (that is, the number of leaves that are directly
attached to the right-most branch).
The product of $\YSym$ does not preserve the
grading by length. For instance,
\[ M_{\epsffile{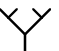}}\cdot M_{\epsffile{figures/1.eps}}=
 M_{\epsffile{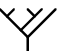}}+ M_{\epsffile{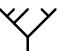}}+ M_{\epsffile{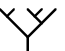}}
 +2\cdot M_{\epsffile{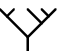}} + M_{\epsffile{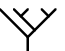}}.\]

Consider the associated graded Hopf algebra, $\gr(\YSym)$. As coalgebras,
$\gr(\YSym)=\YSym$ but the product has been altered by removing terms of lower length
(Section~\ref{S:cofree}). Thus, in $\gr(\YSym)$,
 \[M_{\epsffile{figures/231.eps}}\cdot M_{\epsffile{figures/1.eps}}=
2\cdot M_{\epsffile{figures/3421.eps}}+ M_{\epsffile{figures/4231.eps}}.\]

$\YSym$ admits a Hopf grading, given by the number of internal nodes of a planar
binary tree (one less than the number of leaves). The isomorphism $\YSym\cong Q(V)$
matches this grading with the grading by weight.
This also yields a grading on $\gr(\YSym)$, which corresponds to the grading by weight
under the isomorphism $\gr(\YSym)\cong \Sh(V)$ of Proposition~\ref{P:shuffle-2}.

We relate the graded Hopf algebras $\gr(\YSym)$ and $\GL$ (graded
 by one less than the number of leaves and one less than the number of
nodes, respectively). The dual of $\gr(\YSym)$ is with respect to this grading.

\begin{thm}\label{T:GL-YSym}
 There is an isomorphism of graded Hopf algebras $\Psi:\gr(\YSym)^*\to\GL$ uniquely
determined by
\begin{equation}\label{E:GL-YSym}
M_t^*\mapsto \psi(t)
\end{equation}
for $\iten$-irreducible planar binary trees $t$.
\end{thm}
\begin{proof}
According to the previous discussion,  $\gr(\YSym)$ is the shuffle Hopf algebra on the
subspace $V$ and the number of internal nodes corresponds to the grading by weight.
Therefore $\gr(\YSym)^*$ is the tensor Hopf algebra $T(V^*)$ on the graded
dual space. 
Thus \eqref{E:GL-YSym} determines a  morphism of algebras
$\Psi:\gr(\YSym)^*\to\GL$. Since the number of nodes of $\psi(t)$ is
 the number of leaves of $t$, $\Psi$ preserves the Hopf gradings.
Moreover, $\Psi$ preserves coproducts on a set of algebra generators of
$\gr(\YSym)^*$: the elements $M_t^*$ indexed by $\iten$-irreducible planar binary trees  are primitive generators of the tensor Hopf algebra, and their images
$\psi(t)$ are primitive elements of $\GL$ (since they are planted trees).
Therefore, $\Psi$ is  a morphism of Hopf algebras.
 
We complete the proof by showing that $\Psi$ is invertible.

 Let $t$ be an arbitrary planar binary tree and $t=t_1\iten t_2\iten\cdots\iten t_k$ the
decomposition~\eqref{E:spikes}.
 Then $M^*_t=M^*_{t_1}\cdot M^*_{t_2}\dotsb M^*_{t_n}$, and so
\[
  \Psi(M^*_t)\ =\ \psi(t_1)\cdot\psi(t_2)\dotsb\psi(t_k)\,.
\]
 Since each $t_i$ is planted, Definition~\ref{D:GL} shows that this product is the sum of all
 ordered trees obtained by attaching the root of $\psi(t_{k-1})$ to a node of $\psi(t_k)$, and then
 attaching the root of $\psi(t_{k-2})$ to a node of the resulting tree, and etc.
 The number of children of the root of such a tree is less than $k$, except
 when all the $\psi(t_i)$ are attached to the root, obtaining the ordered tree
 $\psi(t)=\psi(t_1)\iten\psi(t_2)\iten\cdots\iten\psi(t_k)$.

 Linearly ordering both ordered trees and planar binary trees so that
 trees with fewer components precede trees with more components (in the
decompositions~\eqref{E:planted} and~\eqref{E:spikes}), this calculation shows that 
\[\Psi(M^*_t)\ =\ \psi(t)+\text{trees of smaller order.}\] 
Thus $\Psi$ is bijective.
\end{proof}

The main result of Grossman and Larson on the structure of $\GL$~\cite[Theorem 5.1]{GL89} is
contained in the 
proof of Theorem~\ref{T:GL-YSym}. We state it next.
\begin{cor}\label{C:free-ordered}
The set of planted ordered trees freely generates the algebra $\GL$ of ordered trees. Moreover,
$\GL$ is isomorphic to the tensor Hopf algebra on the linear span of the set of planted trees.
\end{cor}
\begin{proof}
As seen in the proof of Theorem~\ref{T:GL-YSym}, $\GL\cong T(V^*)$ as Hopf algebras.
 The isomorphism maps a basis of $V^*$ to the set of planted trees, so the
result follows.
\end{proof}

\begin{rem} We point out that one may construct an isomorphism of graded
Hopf algebras $\gr(\YSym)^*\cong\GL$ from any bijection between the
set of planted trees with $n$ nodes and the set of $\iten$-irreducible planar binary trees with $n$ leaves, instead of the map $\psi$ we used. In fact, 
since they are tensor Hopf algebras, any degree-preserving bijection between the sets of generators determines a unique isomorphism of graded Hopf algebras.

The number of planted trees with $n+2$ nodes (or $\iten$-irreducible planar binary trees with $n+2$ leaves) is the Catalan number $\frac{1}{n+1}\binom{2n}{n}$.
\end{rem}

\section{The Hopf algebra of heap-ordered trees} \label{S:SSym}
We show that the graded dual to $\gr(\SSym)$ is isomorphic to the cocommutative
Hopf algebra of heap-ordered trees defined by Grossman and Larson~\cite{GL89}.

We first review the definition of the Hopf algebra of heap-ordered trees.

A {\em heap-ordered tree} is an ordered tree $x$ together with a labeling of the nodes (a bijection
$\Nod(x)\to\{0,1,\ldots,n\}$) such that:
\begin{itemize}
\item The root of $x$ is labeled by $0$;
\item The labels increase as we move from a node to any of its children;
\item The labels decrease as we move from left to right within the children of each node.
\end{itemize}
The heap-ordered trees with 1, 2, 3, and 4 nodes are shown below:
\[\hzero,\qquad
\hone,\qquad
\htwoone,\htwotwo,\qquad
\hthreeone,\ \hthreetwo,\ \hthreethree,\ \hthreefour,\ \hthreefive,\ \hthreesix.\]

The constructions for ordered trees described in Section~\ref{S:YSym} may be adapted to
the case of heap-ordered trees.

Let $x$ and $y$ be heap-ordered trees. Suppose $x$ has $k$ planted
components (these are ordered trees).
Given a function $f:[k]\to\Nod(y)$, the ordered tree $x\tof y$ may be turned into
a heap-ordered tree by keeping the labels of $y$ and  
incrementing the labels of
 $x$ uniformly by the highest label of $y$.
Given a subset $S=\{i_1<\dotsb<i_p\}\subseteq [k]$, the ordered tree $x_S$
may be turned into a heap-ordered tree by standardizing the labels,
which is to replace the $i$th smallest label by the number $i$, for each $i$.

\begin{defi}\label{D:HGL}
 The Grossman-Larson Hopf algebra $\HGL$ of heap-ordered trees is
 the formal linear span of all heap-ordered trees with product and coproduct as
 follows. Given heap-ordered trees $x$ and $y$ as above, we set
\begin{eqnarray*}
   x\cdot y &=& \sum_{f:[k]\to\Nod(y)}\!\!\!\! x\tof y\,,\\
  \Delta(x) &=& \sum_{S\subseteq[k]}x_S\otimes x_{S^c}\,.\rule{0pt}{16pt}
 \end{eqnarray*}
\end{defi}
For instance, 
\[\htwotwo\cdot\hone\ =\  \hthreetwo+
\hthreefour+\hthreethree+\hthreesix,\]
\begin{align*}
\Delta\bigl(\hfourone\bigr)\ &= \
\hzero\ten\hfourone\ +\
\hone\ten\hthreethree\ +\ \hone\ten\hthreefour\ +\ \htwoone\ten\htwotwo\\
&+\  \
\hfourone\ten\hzero +\ \hthreethree\ten\hone\ +\ \hthreefour\ten\hone\ +\ \htwotwo\ten\htwoone\ .
\end{align*}
$\HGL$ is a graded cocommutative Hopf algebra, where the degree of an ordered tree is one less than
the number of nodes~\cite[Theorem 3.2]{GL89}.

Heap-ordered trees on $n{+}1$ nodes
are in bijection with permutations on $n$ letters. We construct a permutation from such a
tree by listing the labels of all non-root nodes in such way that
the label of a node $i$ is listed to the left of the label of a node $j$ precisely
when $i$ is below or to the left of $j$ (that is, when
 $i$ is a predecessor of $j$, or $i$ is among the left descendants of the
nearest common predecessor between $i$ and $j$).
For instance, the six heap-ordered trees on $4$ nodes above
correspond respectively to $123$, $132$, $213$, $312$, $231$, and $321$.

Let $\phi$ be the inverse bijection. Given a permutation $u$, the heap-ordered tree
$\phi(u)$ is computed  as follows.
Let $u(1),\ldots,u(n)$ be the values of $u$ and set $u(0):=0$.
\begin{itemize}
\item Step 0. Start from a root labeled $0$.
\item Step 1. Draw a child of the root labeled $u(1)$.
\item Step $i$, $i=2,\ldots,n$. Draw a new node labeled $u(i)$. Let
$j\in\{0,\ldots,i{-}1\}$ be the maximum index such that $u(i)>u(j)$. The new node
is a child of the node drawn in step $j$, and it is placed to the right
of any  previous children of that node.
\end{itemize}
For instance,
\[\phi(4231) \ = \ \raisebox{-6pt}{\hfourthree},\quad\mbox{ and }\quad
  \phi(1342) \ = \ \raisebox{-8pt}{\hfourfour}.\]

Given two heap-ordered trees $x$ and $y$, the ordered tree
$x\iten y$ may be turned into a heap-ordered tree by 
incrementing all labels of the nodes in $x$
by the maximum label of a node in $y$. For instance,
\[
 \begin{picture}(50,80)
  \put(0,5){\epsfysize=70pt\epsfbox{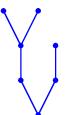}}
   \put( 5,73){\small\Blue{6}} \put(30,73){\small\Blue{4}}  
   \put( 5,45){\small\Blue{3}}    \put(40,48){\small\Blue{5}}
   \put( 5,25){\small\Blue{2}}    \put(40,25){\small\Blue{1}}  
                \put(16, 0){\small\Blue{0}}  
 \end{picture}\ 
  \epsfxsize=0.35in\epsfbox{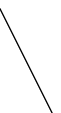}\ 
 \begin{picture}(40,80)
   \put(7,5){\epsfysize=70pt\epsfbox{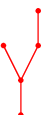}}
                                \put(25,73){\small\Red{4}} 
   \put( 0,48){\small\Red{3}}   \put(23,48){\small\Red{2}}
   \put(12,25){\small\Red{1}} 
   \put(12, 0){\small\Red{0}} 
 \end{picture}
\ \
   \raisebox{24pt}{=}\ 
  \begin{picture}(100,80)
   \put(7,5){\epsfysize=70pt\epsfbox{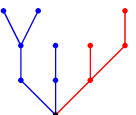}}
   \put(-5,73){\small\Blue{10}}  \put(35,73){\small\Blue{8}}  
   \put(10,45){\small\Blue{7}}   \put(35,48){\small\Blue{9}}
   \put(10,25){\small\Blue{6}}   \put(35,25){\small\Blue{5}}  

                                \put(95,73){\small\Red{4}} 
   \put(70,48){\small\Red{3}}   \put(95,48){\small\Red{2}}
   \put(72,25){\small\Red{1}} 
   \put(35, 0){\small 0} 
 \end{picture}\ .
\]

The operation $\iten$ is  associative on heap-ordered trees, so each such tree has
a unique 
irreducible 
decomposition into $\iten$-irreducible ones. As for ordered trees,
the heap-ordered trees that are planted are $\iten$-irreducible. There are, however,
many other $\iten$-irreducible heap-ordered trees. For instance, while
\[\tone\ \iten\ \ttwoone \ = \ \tthreethree,\]
the heap-ordered tree
\[\hthreethree\]
is $\iten$-irreducible.

The operation $u\iten v$ between permutations~\cite{LR02}  
is obtained by first listing the values of $u$,  
incremented 
by the highest value of $v$, and then listing the
values of $v$ to its right. For instance,
\[\Blue{231}\iten \Red{21}\ = \ \Blue{453}\Red{21}\,. \]
A permutation $w$ has a {\em global descent} at position $p$ if
$w=u\iten v$ with $u$ a permutation of $p$ letters.
Thus, the
$\iten$-irreducible permutations are the permutations with no global descents
(see~\cite[Corollary 6.4]{ASa} for their enumeration). 

The definition of $\phi$ (or its inverse) makes it clear that
\[\phi(u\iten v)=\phi(u)\iten\phi(v)\]
for any permutations $u$ and $v$. 
In particular, $\iten$-irreducible heap-ordered trees correspond
to $\iten$-irreducible permutations under $\phi$.

\smallskip
In~\cite{ASa}, we introduced a linear basis
$M_w$ of $\SSym$, indexed by permutations $w$, which is obtained from the original
basis of Malvenuto and Reutenauer by a process of M\"obius inversion. We showed that $\SSym$ is a
cofree graded  coalgebra and the space $V$ of primitive elements is
the linear span of the elements $M_w$ indexed by 
$\iten$-irreducible 
permutations $w$~\cite[Theorem 6.1, Corollary 6.3]{ASa}.
The isomorphism $Q(V)\cong\SSym$ is
\[M_{w_1}\iten\cdots\iten M_{w_k}\longleftrightarrow M_{w_1\iten\cdots\iten w_k}\,.\]
The resulting grading by length on $\SSym$ is given by the number of $\iten$-irreducible components
in the decomposition of a permutation $w$.
The product of $\SSym$ does not preserve this
grading by length. For instance, in $\SSym$,
\[ M_{231}\cdot M_{1}=
 M_{2314}+M_{2413}+M_{2341}+2\cdot M_{2431}+ M_{3412}+2\cdot M_{3421} + M_{4231}\,.\]
In this product, $M_{231}$ has length $2$, $M_1$ has length $1$, and the only elements of length
$3$ are $M_{3421}$ and $M_{4231}$. Thus, in the associated graded Hopf algebra $\gr(\SSym)$,
 \[M_{231}\cdot M_{1}=2\cdot M_{3421}+M_{4231}\,.\]

$\SSym$ admits a Hopf grading, in which a permutation on $n$ letters has degree $n$.
The isomorphism $\SSym\cong Q(V)$ matches this grading with the grading by weight.
This also yields a grading on $\gr(\SSym)$, which corresponds to the grading by weight
under the isomorphism $\gr(\SSym)\cong \Sh(V)$ of Proposition~\ref{P:shuffle-2}.

Theorem~\ref{T:HGL-SSym} relates the dual of $\gr(\SSym)$ with respect to this grading, with
 $\HGL$, graded by  one less than the number of nodes.

We define the {\it order} of a heap-ordered tree $x$ to be the
pair $(k,l)$, where $k$ is the number of planted components of $x$ and $l$ is
the number of irreducible components of $x$.
We use the following version of the lexicographic order to compare trees:
\[
  (k,l)\ <\ (m,n)\quad\mbox{if}\quad k<m\quad\mbox{or}\quad
   k=m\ \mbox{and}\ l>n\,. 
\]
That is, trees with more planted components have higher order, but among trees
with the same number of planted components, then those with {\it fewer}
irreducible components have higher order. 

Let $x$ be a heap-ordered tree
and $\alpha$ an arbitrary element of $\HGL$. The notation
\[\alpha=x+\tsso\]
indicates that $\alpha-x$ equals a linear combination of heap-ordered trees each of which
is of strictly smaller order than $x$. Not every $\alpha$ can be written in this form, as
 several trees of the same order may appear in $\alpha$.

\begin{lem}\label{L:order} If $\alpha=x+\tsso$ and $\beta=y+\tsso$, then
$\alpha\cdot\beta=x\iten y+\tsso$
\end{lem}
\begin{proof} Consider first the product of two heap-ordered trees $x'$ and $y'$
having orders $(k,l)$ and $(m,n)$  respectively.
This is the sum of all heap-ordered trees
 obtained by attaching the  planted components 
 of $x'$ to nodes of $y'$. 
 Every such tree will have fewer than $k+m$ planted components, except the
 tree obtained by attaching all planted components of $x'$ to the root of
 $y'$, which will be $x'\iten y'$ and will have $k+m$ planted components.

Therefore, among the trees appearing in $\alpha\cdot\beta$, the ones with the maximum number of
planted components are those of the form $x'\iten y'$, with $x'$ and $y'$ having the same numbers
of planted components as $x$ and $y$ respectively. 
Among these we find the tree $x\iten y$.
For any of the remaining trees with the maximum number of
planted components, either $x'$ has more irreducible components than $x$ or $y'$
has more irreducible components than $y$, by hypothesis. Since the number of irreducible
components of $x'\iten y'$ is $l+n$, the tree $x'\iten y'$ has more irreducible components than
$x\iten y$, and hence it is of smaller order.
\end{proof}

Applying Lemma~\ref{L:order} inductively we deduce that any heap-ordered tree $x$ is the leading
term in the product of its irreducible components.
This implies that the the set of
irreducible heap-ordered trees freely generates the algebra 
$\HGL$ of ordered trees. 
This result is due to Grossman and Larson~\cite[Theorem 6.3]{GL89}.
Irreducible heap-ordered trees are not necessarily primitive.
We refine this result of Grossman and Larson, giving primitive generators and
relating the structure of $\HGL$ explicitly to that of $\gr(\SSym)^*$.

\medskip

We assume from now on that the base field $\field$ is of characteristic $0$.

\medskip 

We need one more tool: the {\em first Eulerian idempotent}~\cite{GS91},
~\cite[Section 4.5.2]{Lod98},~\cite[Section 8.4]{Re93}. For any graded connected Hopf algebra $H$, the identity map $\id:H\to H$
is locally unipotent with respect to the convolution product of $\End(H)$. 
Here $1$ denotes the composite $H\map{\epsilon}\field\map{u}H$ of the counit and unit maps of $H$ (the unit element for the convolution product). Therefore,
\[\euler:=\log(\id)=\sum_{n\geq 1}\frac{(-1)^{n+1}}{n}(\id-1)^{\ast n}\]
is a well-defined linear endomorphism of $H$. The crucial fact is that if $H$ is
cocommutative, this operator is a projection onto the space of primitive elements of
$H$: $\euler:H\onto P(H)$~\cite{Pat},~\cite[pages 314-318]{Sch94}.

\begin{lem}\label{L:euler} Let $x$ be a $\iten$-irreducible heap-ordered tree. Then
\[\euler(x)=x+\tsso\]
\end{lem}
\begin{proof} 
In any graded connected Hopf algebra $H$, the map $\id -1$ is the projection of $H$ onto the part 
of positive degree, and the convolution power $(\id-1)^{\ast n}$ equals the map 
$m^{(n-1)}\circ(\id-1)^{\otimes n}\circ\Delta^{(n-1)}$. 
Let $x$ be a heap-ordered tree with $k$ planted components.
Iterating the coproduct of $\HGL$ (Definition~\ref{D:HGL})
 gives
 \[\Delta^{(n-1)}(x)\ =\ \sum_{S_1\sqcup\dotsb\sqcup S_n=[k]}
 x_{S_1}\otimes\dotsb\otimes x_{S_n}\,,\]
the sum over all ordered decompositions of $[k]$ into $n$ disjoint subsets. Applying
$(\id-1)^{\otimes n}$ to this sum has the effect of erasing all terms corresponding to 
decompositions involving at least one empty set. Therefore,
 \[(\id -1)^{\ast n}(x)=\sumsub{S_1\sqcup\cdots\sqcup S_n= [k]\\S_i\neq\emptyset}
 x_{S_1}\dotsb x_{S_n}\,,\]
the sum now over all {\em set-compositions} of $[k]$ (decompositions into non-empty disjoint
subsets). In particular, this sum is $0$  when $n>k$. Thus,
\[\euler(x)\ =\ \sum_{n=1}^k\frac{(-1)^{n+1}}{n}\sumsub{S_1\sqcup\cdots\sqcup S_n=
[k]\\S_i\neq\emptyset}
 x_{S_1}\dotsb x_{S_n}\,.\]
By Lemma~\ref{L:order}, $x_{S_1}\dotsb x_{S_n}=x_{S_1}\iten\dotsb\iten x_{S_n}+\tsso$
Each tree $x_{S_1}\iten\dotsb\iten x_{S_n}$ has $k$ planted components (as many as $x$)
and at least $n$ irreducible components. Hence, among these trees, the one of highest order
 is $x$, which corresponds to the trivial decomposition of $[k]$ into $n=1$ subset.
Thus, among all trees appearing in $\euler(x)$, there is one of
highest order and it is $x$. 
\end{proof}

For example, if $x=\raisebox{-3pt}{\hthreethree}$, then
\[\euler(x)=\raisebox{-3pt}{\hthreethree}
-\frac{1}{2}\Bigl(\raisebox{-3pt}{\ \hthreefour}+\raisebox{-3pt}{\hthreefive}+
\raisebox{-3pt}{\hthreetwo}\Bigl)-
\raisebox{-3pt}{\hthreeone}\,.\]
The tree $x$ is of order $(2,1)$, the next two trees are of order $(2,2)$, and the last two of order
$(1,1)$.

 \begin{thm}\label{T:HGL-SSym} Assume $\ch(\field)=0$.
 There is an isomorphism of graded Hopf algebras $\Phi:\gr(\SSym)^*\to\HGL$ uniquely
determined by
\begin{equation}\label{E:HGL-SSym}
M_w^*\longmapsto \euler\bigl(\phi(w)\bigr)
\end{equation}
for $w$ a $\iten$-irreducible permutation.
\end{thm}
\begin{proof}
By 
the discussion preceding Lemma~\ref{L:order}, $\gr(\SSym)^*\cong T(V^*)$. Therefore,~\eqref{E:HGL-SSym}
determines a morphism of graded algebras $\Phi$. Since $\HGL$ is cocommutative, $\euler\bigl(\phi(w)\bigr)$ is a primitive
element of $\HGL$. Thus $\Phi$
preserves primitive elements and hence it is a morphism of Hopf algebras.

It remains to verify that $\Phi$ is invertible.  Let $w=w_1\iten\dotsb\iten w_k$ be the
irreducible decomposition of a permutation $w$. 
 Let $x:=\phi(w)$. Since $\phi$ preserves the operations $\iten$,
 the irreducible components of $x$ are $x_i:=\phi(w_i)$,
  $i=1,\ldots,k$.
 On the other hand, $M^*_w=M^*_{w_1}\dotsb M^*_{w_k}$, so
\[\Phi(M^*_w)\ =\ \euler\bigl(x_1\bigr)\dotsb\euler\bigl(x_k\bigr)\,.\]
{}From Lemmas~\ref{L:order} and~\ref{L:euler} we deduce
\[
\Phi(M^*_w)=x_{1}\iten\dotsb\iten x_{k}+\tsso=\phi(w)+\tsso
\]
 As in the proof of Theorem~\ref{T:GL-YSym},
 this shows that $\Phi$ is invertible, by triangularity.
\end{proof}

 Let $W$ be the graded space where $W_n$ is spanned by the elements $M^*_w$, for $w$ an
irreducible permutation of $[n]$.  
 {}From the proof of Theorem~\ref{T:HGL-SSym}, we deduce the following corollary.

\begin{cor}\label{C:free-heap}
 The Hopf algebra $\HGL$ of ordered trees is isomorphic to the tensor Hopf algebra on
 a graded space  $W=\oplus_{n\geq 0}W_n$ with $\dim W_n$ equal to the number of 
 irreducible heap-ordered trees on $n+1$ nodes  
 (or the number of irreducible permutations of $[n]$).
\end{cor}

\begin{rem} As pointed out to us by Holtkamp, the  use of the Eulerian idempotent 
in Theorem~\ref{T:HGL-SSym}
is similar to that
encountered in a proof of the Milnor-Moore theorem~\cite[Theorem 5.18]{MM}, ~\cite[Theorem 4.5]{Qui}.
\end{rem}



\end{document}